\documentclass[11pt,oneside]{amsart}
\usepackage{amscd,amssymb}

\usepackage{graphicx}
\theoremstyle{plain}
\newtheorem{Thm}[equation]{Theorem}
\newtheorem{Ex}[equation]{Example}

\newtheorem{Rmk}[equation]{Remark}

\newtheorem{Prop}[equation]{Proposition}

\newtheorem{Lem}[equation]{Lemma}
\newtheorem{Def}[equation]{Definition}
\numberwithin{equation}{section}

\newcommand{\e}{\epsilon}

\newcommand{\br}{\mathbb{R}}

\newcommand{\ba}{\backslash}

\newcommand{\G}{\Gamma}

\renewcommand{\P}{\mathcal P}

\newcommand{\la}{\langle}
\newcommand{\ra}{\rangle}

\newcommand{\bp}{\begin{pmatrix}}
\newcommand{\ep}{\end{pmatrix}}

\newcommand{\SO}{\operatorname{SO}}

\newcommand{\BMS}{\op{BMS}}
\newcommand{\BR}{\op{BR}}
\newcommand{\bi}{\begin{itemize}}

\newcommand{\ei}{\end{itemize}}

\newcommand{\B}{\mathbb B}

\newcommand{\op}{\operatorname}

\newcommand{\PS}{\operatorname{PS}}
\newcommand{\Leb}{\operatorname{Leb}}

\newcommand{\T}{\op{T}}
\renewcommand{\PS}{\operatorname{PS}}
\renewcommand{\BR}{\operatorname{BR}}
\renewcommand{\Leb}{\operatorname{Leb}}
\renewcommand{\BMS}{\operatorname{BMS}}

\newcommand{\bB}{\mathbb B}
\renewcommand{\S}{\mathbb S^2}
\newcommand{\W}{\mathcal W}

\newcommand{\N}{\mathcal N}
\newcommand{\ee}{\end{equation}}\newcommand{\be}{\begin{equation}}
\newcommand{\ees}{\end{equation*}}\newcommand{\bes}{\begin{equation*}}

\begin{document}

\title[Counting  visible circles on the sphere]
{Counting  visible circles on the sphere and Kleinian groups}
%{\tiny\it\rm draft}

%\title[Asymptotic distribution of circles on the sphere]
%{The asymptotic distribution of circles in the orbits of Kleinian groups
% on the sphere}
%{\tiny\it\rm draft}

\author{Hee Oh and Nimish Shah}

\address{Mathematics department, Brown university, Providence, RI, USA
and Korea Institute for Advanced Study, Seoul, Korea}
\email{heeoh@math.brown.edu}
\address{The Ohio State University, Columbus, OH and Tata Institute of Fundamental Research, Mumbai, India}
\email{shah@math.osu.edu}

\thanks{Oh is partially supported by NSF
   grant 0629322}
%\email{heeoh@math.brown.edu}
%\address{Mathematics department, Brown university, Providence, RI
%and Korea Institute for Advanced Study, Seoul, Korea}
\begin{abstract}
For a circle packing $\P$ on the sphere invariant
under a geometrically finite Kleinian group, we compute the asymptotic
of the number of circles in $\P$ of spherical curvature at most $T$
which are contained in any given region.
\end{abstract}

%\email{heeoh@math.brown.edu}
%\address{Mathematics department, Brown university, Providence, RI
%and Korea Institute for Advanced Study, Seoul, Korea}

\maketitle
%\tableofcontents

\section{Introduction}
In the unit sphere $\mathbb S^2=\{x^2+y^2+z^2=1\}$
with the Riemannian metric
induced from $\br^3$, the distance between two points is simply the angle between the rays connecting
them to the origin $o$.
Let $\P$ be a circle packing on the sphere ${\mathbb S^2}$, i.e., a union of circles
which may intersect with each other.

In the beautiful book {\it Indra's pearls}, Mumford, Series and Wright ask
the question (see \cite[5.4 in P. 155]{MumfordSeriesWright})
\begin{equation*}\label{q} \text{\it How many visible circles are there?}\end{equation*}
%in the limit set of a finitely generated Schottky group of $\PSL_2(\c)=\op{Isom}_+(\mathbb B)$.
%As  $\PSL_2(\c)$

The {\it visual size} of a circle $C$ in $\mathbb S^2$ can be measured by its
spherical radius $0<\theta(C)\le \pi/2$, that is, the half of the visual angle of $C$
from the origin $o=(0,0,0)$. We label the circles by their spherical curvatures given by
$$\op{Curv}_{S}(C):=\cot \theta(C).$$
%noting that $\op{Curv}_{S}(C)\to \infty$ as $\theta(C)\to 0$.

We suppose that $\P$ is infinite and locally finite in the sense that for any $T>1$,
 there
are only finitely many circles in $\P$ of spherical curvature at most $T$.
We then set for any subset $E\subset \mathbb S^2$,
$$N_T({\mathcal P}, E)=\{C\in \mathcal P: C\cap E\ne\emptyset,\;\;
 \op{Curv}_{S}(C) <T \} <\infty .$$

In order to present our main result on the asymptotic for $N_T(\P, E)$,
we consider the Poincare ball model $\mathbb B=\{x_1^2+x_2^2+x_3^2<1\}$
 of the hyperbolic $3$-space with the metric given by
$\frac{2\sqrt{dx_1^2+dx_2^2+dx_3^2}}{1-(x_1^2+x_2^2+x_3^2)}.$
The geometric boundary of $\mathbb B$ naturally identifies with ${\mathbb S^2}$.
% and that
%for any circle $C$ in $\mathbb S^2$,
%we have
%\begin{equation}\label{eq} \sin \theta(C)=\frac{1}{\cosh d(\hat C, o)} \end{equation}
% where $\hat C\subset \mathbb B$ is the convex hull of $C$.

%The convex core $C_\G\subset \G\ba \mathbb B$ of $\G$ is defined to be the minimal convex
%set in $\mathbb B$ mod $\G$ which contains all geodesics connecting
%any two points in $\Lambda(\G)$.
%A Kleinian group $\G$ is geometrically finite
%if and only
% if the unit neighborhood of its convex core has finite volume.
% $\G$ is called {\it convex co-compact}
% if its convex core is compact.

\begin{figure}
 \begin{center}
    \includegraphics[height=5cm]{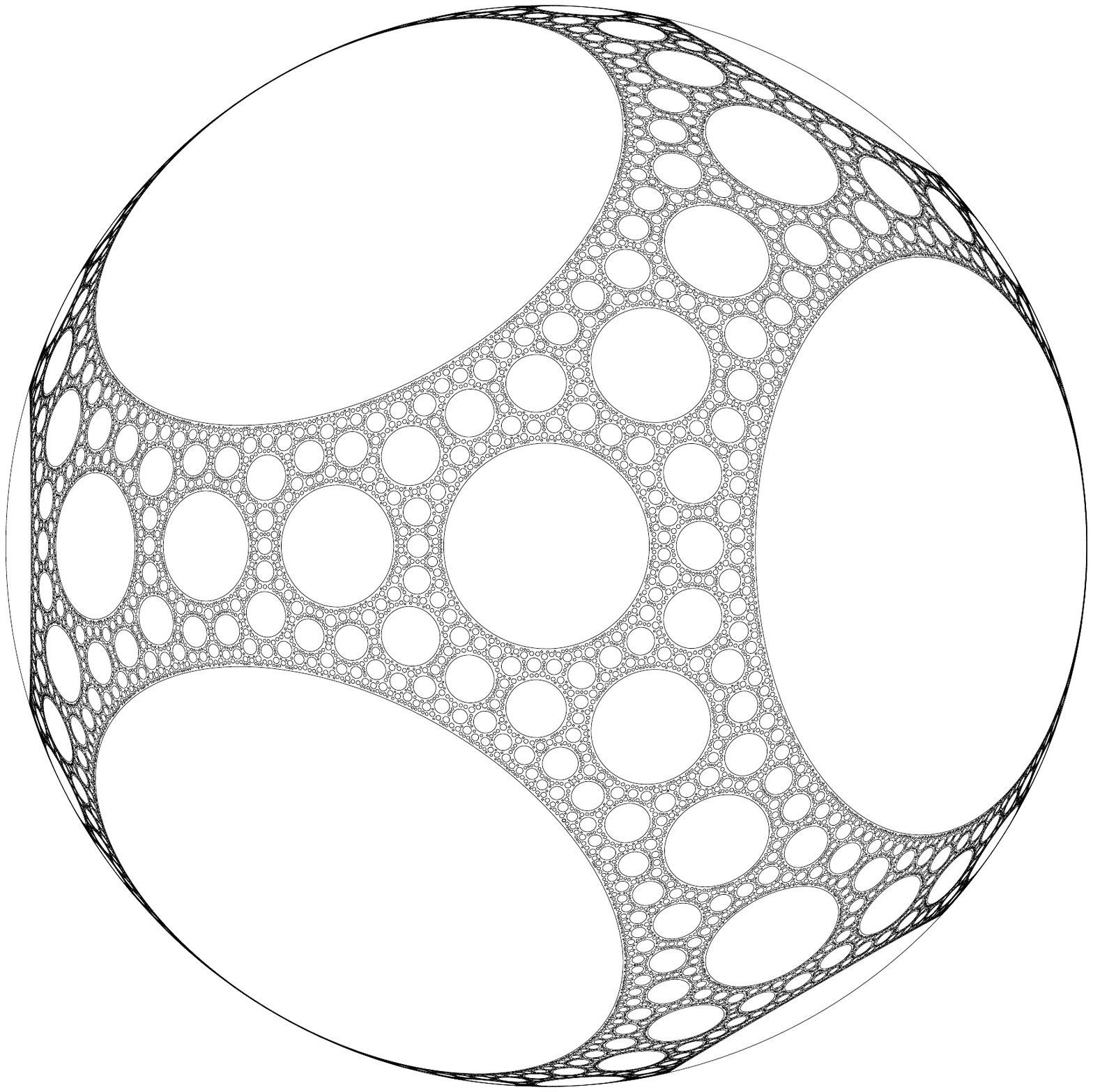}
 \includegraphics[height=5cm]{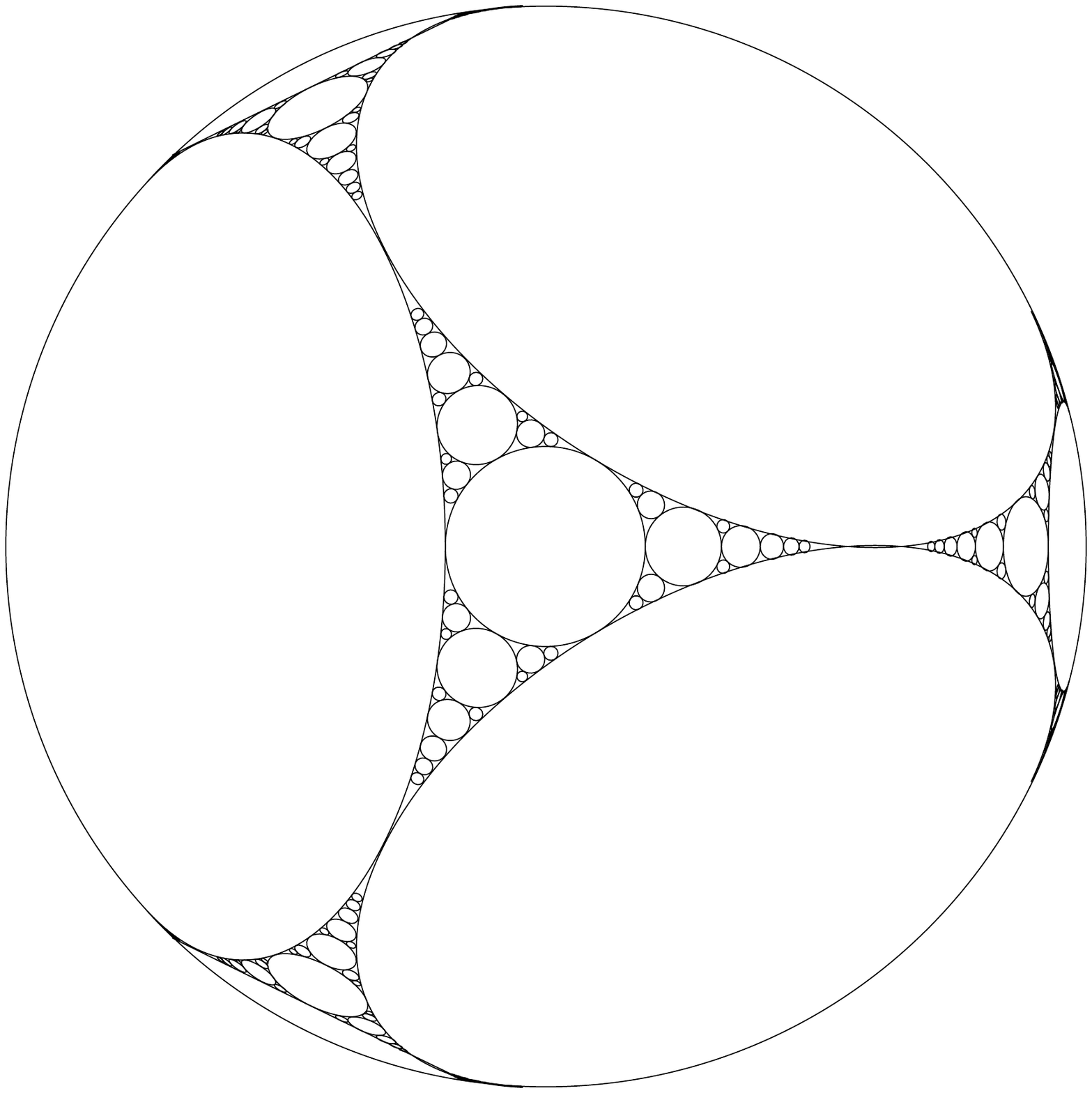}
     \caption{Sierpinski curve and Apollonian gasket (by C. McMullen)}
    \label{f1}
 \end{center}
\end{figure}

Let $G$ denote the group of orientation preserving isometries of
 $\mathbb B$. A torsion-free discrete subgroup $\G$ of $G$ is called a Kleinian group.
A Kleinian group $\G$ is called {\it geometrically finite} if
$\G$ admits a finite sided fundamental domain in $\mathbb B$, and
{\it non-elementary} if $\G$ has no abelian subgroup of finite index.
We denote by $\Lambda(\G)\subset\mathbb S^2$ the limit set of $\G$,
 that is, the set
of accumulation points of an orbit of $\G$ in $\mathbb B\cup\mathbb S^2$.
The critical exponent $0\le \delta_\G \le 2$ of $\G$ is known to be
 positive for $\G$ non-elementary and equal to
 the Hausdorff dimension of the limit set $\Lambda(\G)$ for $\G$ geometrically finite \cite{Sullivan1979}.

%Denote by $\delta_\G$ the critical exponent of $\G$: the abscissa of convergence for
%the Poincare series $$\mathcal P_\G(s)=\sum_{\gamma\in \G} e^{-s d(x,\gamma x)}$$
%for $x\in \mathbb B$. We note that for $\G$ non-elementary, $\delta_\G $ is positive
%and for $\G$ geometrically finite, $\delta_\G$ is equal to the Hausdorff dimension of $\Lambda(\G)$
%\cite{Sullivan1979}.
%In our setting, $1\le \delta_\G <2$, as $\delta_\G=2$ implies $\Omega_\G=\emptyset$ by Sullivan.
%We assume in the following that $\G$ is geometrically finite and non-elementary, i.e.,
 %$\G$ has no abelian subgroup of finite index.
 % Generalizing the work of Patterson  \cite{Patterson1976}
%for the Poincare disk, Sullivan \cite{Sullivan1979} constructed a
%$\G$-invariant conformal density $\{\nu_x:x\in \mathbb B\}$ on $\Lambda(\G)$ of
%dimension $\delta_\G$ which is unique up to homothety.

For a vector $u$ in the unit tangent bundle $\T^1(\B)$, denote by $u^{+}\in \S$
  the forward end point of the geodesic determined by $u$, and by $\pi(u)\in \B$ the basepoint of $u$.
 For $x_1,x_2\in \B$ and $\xi\in \S$, $\beta_\xi(x_1,x_2)$ denotes the signed distance between horospheres
 based at $\xi$ and passing through $x_1$ and $x_2$.

For a non-elementary geometrically finite Kleinian group
$\G$, denote by $\{\nu_x:x\in \mathbb B\}$ the Patterson-Sullivan density for $\G$ (\cite{Patterson1976}, \cite{Sullivan1979}), i.e., a $\G$-invariant conformal density  of
dimension $\delta_\G$ on $\Lambda(\G)$,  which is unique up to homothety.
 \begin{Def}[The $\G$-skinning size of $\P$]
 {\rm For a circle packing $\P$ on $\S$ invariant under $\G$, we define $0\le \op{sk}_\G(\P)\le\infty$ as follows:
$$\op{sk}_\G(\P):=\sum_{i\in I} \int_{s\in \op{Stab}_{\G} (C_i^\dagger)\ba C_i^\dagger}  e^{\delta_\G
\beta_{s^+}(x,\pi(s))}d\nu_{x}(s^+)$$
where $x\in \mathbb B$, $\{C_i:i\in I\}$ is a set of representatives of $\G$-orbits in $\P$ and
$C_i^\dagger\subset \T^1(\mathbb B)$ is the set of unit normal vectors to the convex hull of $C_i$. }\end{Def}
By the conformal property of $\{\nu_x\}$, the definition of
 $\op{sk}_\G(\P)$ is independent of the choice of $x$ and the choice of representatives $\{C_i\}$.

%The set $\Omega(\G):=\mathbb S^2-\Lambda(\G )$ is called the domain of discontinuity for $\G$.

\begin{Thm}\label{m1}
Let $\G$ be a non-elementary geometrically finite Kleinian group and
 $\mathcal P=\cup_{i\in I} \Gamma(C_i)$ be an infinite,
  locally finite, and  $\G$-invariant circle packing on the sphere ${\mathbb S^2}$ with finitely many $\G$-orbits.

Suppose one of the following conditions hold:
\begin{enumerate}
\item $\G$ is convex co-compact{\footnote{A discrete subgroup $\G<G$ is called convex co-compact
if the convex core $C_\G$ (that is, the minimal convex
set in $\G\ba \bB$ containing all geodesics connecting
any two points in $\Lambda(\G)$)
  is compact.}};
\item  all circles in $\P$ are mutually disjoint;
\item $\cup_{i\in I} C_i^\circ \subset \Omega(\G)$ where $C_i^\circ$ denotes the open disk enclosed by $C_i$ and $\Omega(\G)=\S-\Lambda(\G )$ is  the domain of discontinuity for $\G$.
\end{enumerate}
Then for any Borel subset $E\subset \mathbb S^2$ whose boundary has zero Patterson-Sullivan measure,
 $$N_T(\P, E)
\sim
 \frac{ \op{sk}_\G(\P)}{\delta_\G \cdot |m^{\BMS}_\G|} \cdot \nu_o(E)\cdot
 (2 T)^{\delta_\G}  \quad\text{as $T\to \infty$}$$
 where $o=(0,0,0)$, $0<\op{sk}_\G(\P)<\infty$ and
  $0<|m^{\BMS}_\G|<\infty$ is the total mass
of the Bowen-Margulis-Sullivan measure associated to $\{\nu_x\}$ (Def. \ref{BMS}) on $\T^1(\G\ba \B)$.
\end{Thm}

\begin{Rmk}\rm
The Patterson-Sullivan density is known to be atom-free, and hence
the above theorem works for any Borel subset $E$
intersecting $\Lambda(\G)$ in finitely many points. If $\G$ is Zariski dense in $G$, then
any proper real subvariety of $\mathbb S^2$ has zero Patterson-Sullivan density \cite[Cor. 1.4]{FlaminioSpatzier} and hence
Theorem~\ref{m1}
holds for any Borel subset of $\mathbb S^2$ whose boundary is contained in a countable union of real algebraic curves.
\end{Rmk}

\begin{Ex}
{\rm \begin{enumerate}
\item If $X$ is a
finite volume hyperbolic $3$ manifold with
  totally geodesic boundary,
 its fundamental group $\G:=\pi_1(X)$ is geometrically finite and $X$ is homeomorphic
to $\G\ba \mathbb B\cup \Omega(\G)$ \cite{Kojima1992}.
The universal cover $\tilde X$ developed in $\mathbb B$
has geodesic boundary components which are Euclidean hemispheres
normal to $\mathbb S^2$.
Then $\Omega(\G)$  is the union of
a countably many disjoint open disks corresponding to
the geodesic boundary components of $\tilde X$.
The Ahlfors finiteness theorem \cite{Ahlfors1964}
%$\G\ba \Omega(\G)$ is of finite type. This implies
%that
implies that the circle packing $\P$ on $\mathbb S^2$ consisting of
the geodesic boundary components of $\tilde X$  is locally finite  and
has finitely many $\G$-orbits.
Hence provided $\P$ is infinite,
our theorem \ref{m1} applies to counting circles in $\P$.

\item Starting with four mutually tangent circles on the sphere $\S$,
one can inscribe into each of the curvilinear triangle
a unique circle by an old theorem of Apollonius of Perga (c. BC 200).
Continuing to inscribe the circles this way, one obtains an Apollonian circle packing on $\S$ (see Fig. \ref{f1}).
Apollonian circle packings are examples of circle packing obtained in the way described in (1)
(cf. \cite{ErickssonLagarias2007} and \cite{LagariasMallowsWilks}.).
  In the case when $\pi_1(X)$ is convex co-compact, then
no disks in $\Omega(\G)$ are tangent to each other and $\Lambda(\G)$ is known to
be homeomorphic to a Sierpinski curve \cite{Claytor1934} (see Fig. \ref{f1}).

\item
Take $k\ge 1$ pairs of mutually disjoint closed disks
$\{(D_i, D_i') : 1\le i\le k\}$ in $\mathbb S^2$ and choose
 $\gamma_i\in G$ which maps $D_i$ and $D_i'$ and sends the interior of
$D_i$ to the exterior of $D_i'$. The group, say, $\G$,
 generated by $\{\gamma_i\}$ is called a Schottky group of genus $k$ (cf. \cite[Sec. 2.7]{MardenOutercircles}).
The $\G$-orbit of the disks nest down onto the limit set $\Lambda(\G)$ which is totally disconnected.
If we set $\P:=\cup_{1\le i\le k} \G(C_i)\cup \G(C_i')$
where $C_i$ and $C_i'$ are the boundaries of $D_i$ and $D_i'$ respectively,
$\P$ is locally finite, as the nesting disks will become smaller and smaller (cf. \cite[4.5]{MumfordSeriesWright}).
The common exterior of hemispheres above the initial disks $D_i$ and $D_i'$
is a fundamental domain for $\G$ in $\mathbb B$ and hence
$\G$ is geometrically finite.
Since $\P$ consists of disjoint circles,
Theorem \ref{m1} applies to counting circles in $\P$, called Schottky dance
(see \cite[Fig. 4.11]{MumfordSeriesWright}).
%\item When the boundary of $X$ is not entirely totally geodesic,
%$\Omega(\G)$ need not consist only of open disks (see pictures in P. 246 of \cite{MumfordSeriesWright}).
\end{enumerate}}
\end{Ex}

 Since for $o=(0,0,0)$,
 $\sin \theta(C)=\frac{1}{\cosh d(\hat C, o)} $
for the convex hull $\hat C$ of $C$ (cf. \cite[P.24]{Thurstonbook}),
we deduce
$$\op{Curv}_S(C)=\sinh  d(\hat C, o).$$

% As both $\sin \theta $ and $\cosh d$ are monotone functions for $0\le \theta\le \pi/2$
% and $d \ge 0$ respectively,
Hence Theorem \ref{m1} follows from the following:
\begin{Thm}\label{t_main}
Keeping the same assumption as in Theorem \ref{m1},
we have, for any $o\in \mathbb B$,
$$\#\{C\in \P: C\cap E\ne \emptyset, \; d(\hat C, o)<t\}\sim \frac{ \op{sk}_\G(\P)}{\delta_\G \cdot  |m^{\op{BMS}}_\G|}
\cdot \nu_o(E)\cdot e^{\delta_\G \cdot t} \quad\text{as $t\to \infty$} .$$
\end{Thm}

The main result in this paper was announced in \cite{OhICM}
and an analogous problem of counting circles
in a circle packing of the {\it plane} was studied in \cite{KontorovichOh} and \cite{OhShahcircle}.

%A circle counting problem for the bounded Apollonian circle packings in the plane
%was studied in \cite{KontorovichOh} where a similar asymptotic was obtained. In
%a subsequent paper \cite{OhShahcircle}, we consider the counting problem
%for general circle packings in the plane.

\bigskip
\noindent{\bf Acknowledgment:}
We are very grateful to Curt McMullen for generously sharing his intuition and ideas.
 The applicability of our
other paper \cite{OhShahGFH} in the question addressed in this paper came
 up in the conversation
of the first named author with him.
We also thank Yves Benoist, Jeff Brock and Rich Schwartz for useful conversations.
\bigskip

\section{Preliminaries and expansion of a hyperbolic surface}\label{not}
%We present the proof of Theorem \ref{t_main}. Without loss of generality, we may assume that
 %$\mathcal P=\G (C_0)$ for some circle $C_0$ in $\mathbb S^2$.
In this section, we set up notations
as well as recall a result from \cite{OhShahGFH} on the
weighted equidistribution of
expansions of a hyperbolic surface by the geodesic flow.

Denote by $G$  the group of orientation preserving isometries of $\mathbb B$ and
fix a circle $C_0\subset \mathbb S^2$.
Denote by $\hat C_0\subset \B$ the convex hull of $C_0$.
Fix $p_0\in \hat{C_0}$ and $o\in \mathbb B$.
 As $G$ acts transitively on $\mathbb B$, there exists
 $g_0\in G$ such that $$o=g_0(p_0).$$
Denote by $K$ the stabilizer subgroup of $p_0$ in $G$ and by $H$ the stabilizer subgroup
of $\hat C_0$ in $G$.
We note that $H$ is locally isomorphic to $\SO(2,1)$ and has two connected components
one of which is the orientation preserving isometries of $\hat C_0$.
There exist commuting involutions $\sigma$ and $\theta$ of $G$
 such that the Lie subalgebras $\frak h=\op{Lie}(H)$ and $\frak k=\op{Lie}(K)$
are the $+1$ eigenspaces of $d\sigma$ and $d\theta$ respectively.
With respect to the symmetric bilinear form on $\frak g=\op{Lie}(G)$ given by
$$B_\theta(X,Y)=\op{Tr}(ad(X)\circ ad(\theta(Y)),$$
we have the orthogonal decomposition
$$\frak g=\frak k\oplus \frak p =\frak h \oplus \frak q$$
where $\frak p$ and $\frak q$ are
the $-1$ eigenspaces of $d\sigma$ and $d\theta$ respectively.
Let $\frak a$ be a one dimensional subalgebra
of $\frak p \cap \frak q$, $A:=\exp(\frak a)$, and $M$ the centralizer of $A$ in $K$.
The map $K\times \frak p\to G$ given by $(k, X)\mapsto k \exp X$ is a diffeomorphism and for the canonical
projection $\pi: G\to G/K=\mathbb B$, the differential $d\pi:\frak p\to \T_{p_0}(G/K)=\T_{p_0}(\mathbb B)$ is an isomorphism.

%We set $\frak a^+=\exp (\br_{\ge 0}X_0)$
%and $A^+:=\exp \frak a^+$.

%so that $d\pi(X_0)$ is a unit tangent vector
%at $p_0$ toward the interior of the manifold $\tilde X$.
Choosing an element $X_0\in \frak a$ of norm one,
we can identify the unit tangent bundle $\op{T^1}(\mathbb B)$ with $G.X_0=G/M$:
here $g. X_0$ is given by $d\lambda(g) (X_0)$ where $\lambda(g): G\to G$ is the left translation $\lambda(g)(g')=
gg'$ and $d\lambda$ is its derivative at $p_0$.

Setting $A^+=\{ \exp(t X_0): t\ge 0\}$ and
$A^-=\{ \exp(t X_0): t\le 0\}$,
 we have the following generalized Cartan decompositions (cf. \cite{Schlichtkrull1984}):
$$G=KA^-K= HA^+K .$$
in the sense that every element of $g\in G$ can be written as
$g=k_1 a_s k_2= ha_t k$, $s\le 0$, $t\ge 0$, $h\in H,, k_1, k_2, k\in K$.
Moreover, $k_1a_sk_2= k_1'a_{s'}k_2'$ implies
 $s=s'$, $k_1 = k_1'm_1$, and $k_2= m^{-1}_1 k_2'$ for some $m_1\in M$, and
 $ha_tk=h'a_{t'}k'$ implies that $t=t'$, $h = h'm_2$, and $k= m^{-1}_2 k'$ for some
$ m_2\in H\cap K$.

The set $K.X_0=K/M$ represents the set of unit tangent vectors at $p_0$, and
 as $X_0$ is orthogonal to $\frak h\cap \frak p=\T_{p_0}(\hat C_0)$,
 $H.X_0=H/M$ corresponds to the set of unit normal vectors
   to the convex hull $\hat C_0=H/H\cap K$, which will be denoted by $C_0^\dag$.
Moreover if $a_t=\exp(tX_0)$, the set $(H/M)a_t=(Ha_tM)/M$
 represents the orthogonal translate of $\hat C_0$ by distance $|t|$.
We refer to \cite{OhShahGFH} for the above discussion.

Let $\G<G$ be a non-elementary geometrically finite discrete subgroup of $G$ in the rest of this section.
 Recall that $g\in G$ is parabolic
if and only if $g$ has a unique fixed point in $\mathbb S^2$.
\begin{Prop}[\cite{OhShahcircle}] \label{all}\begin{enumerate}
\item
If $\G(C_0)$ is infinite, then $[\G:H\cap \G]=\infty$.

\item  $\G(C_0)$ is locally finite if and only if
the natural projection map $\G\cap H \ba \hat C_0 \to \G\ba \mathbb B$ is proper.

\item  Suppose one of the following conditions hold:
\begin{enumerate}

\item $\G$ is convex co-compact;
\item $\gamma(C_0)$'s are disjoint for $\gamma\in \G$;
\item $C_0^{\circ}\subset \Omega(\G)$ for the open disk $C_0^\circ$ enclosed by $C_0$.
\end{enumerate}
Then  any parabolic element of $\Gamma$
  fixing a point on $C_0$ stabilizes $C_0$.
%$$\Lambda_p(\G)\cap C\subset \Lambda_p(H\cap \G) .$$
\end{enumerate} \end{Prop}

We denote by $\{\nu_x: x\in \mathbb B\}$ the Patterson-Sullivan density for $\G$,
which is unique up to homothety: for any $x,y\in \bB$, $\xi \in \mathbb S^2$ and $\gamma\in
\G$,
$$\gamma_*\nu_x=\nu_{\gamma x};\quad\text{and}\quad
 \frac{d\nu_y}{d\nu_x}(\xi)=e^{-\delta_\G \beta_{\xi} (y,x)}, $$
where $\gamma_*\nu_x(R)=\nu_x(\gamma^{-1}(R))$
and the Busemann function $\beta_\xi(y_1, y_2)$ is given by
$\lim_{t\to\infty} d(y_1, \xi_t)-d(y_2,\xi_t)$ for a geodesic ray
$\xi_t$ toward  $\xi$.

 %We denote by $\{g^t\}$ the geodesic flow on the unit tangent bundle $\T^1(\mathbb B)$
 For $u\in \T^1(\mathbb B)$, we define $u^+\in \mathbb S^2$ (resp. $u^-\in \S$)
 the forward (resp. backward) endpoint of the geodesic determined by $u$ and $\pi(u)\in \bB$
 the basepoint.
 The map $$u\mapsto (u^+, u^-,
\beta_{u^-}(\pi(u), o))$$ yields a homeomorphism between
$\T^1(\mathbb B)$ with $(\mathbb S^{2}\times \mathbb S^{2} - \{(\xi,\xi):\xi\in
\mathbb S^{2}\}) \times \br .$
%Recall that $\{\nu_x:x\in \B\}$ denotes a Patterson-Sullivan density of $\G$.
\begin{Def}[The Bowen-Margulis-Sullivan measure]\label{BMS}
\rm The Bowen-Margulis-Sullivan measure $m^{\BMS}_\G$ (\cite{Bowen1971}, \cite{Margulisthesis}, \cite{Sullivan1984}) associated to $\{\nu_x\}$ is the measure on $\T^1(\G\ba \mathbb B)$ induced by the following
$\G$-invariant measure on $\T^1(\mathbb B)$: for $x\in \mathbb B$,
$$d  \tilde m^{\BMS}(u)=e^{\delta_\G \beta_{u^+}(x, \pi(u))}\;
 e^{\delta_\G \beta_{u^-}(x,\pi(u)) }\;
d\nu_x(u^+) d\nu_x(u^-) dt .$$
\end{Def}
 It follows from the conformality of $\{\nu_x\}$ that
this definition is independent of the choice of $x$.
Sullivan showed that $m^{\BMS}_\G$ is ergodic for the geodesic flow and that the total mass
$|m^{\BMS}_\G|$ is finite \cite{Sullivan1984}.

We denote by $\{m_x:x\in \mathbb B\}$ a
$G$-invariant conformal
density of $\mathbb S^2$ of dimension $2$, which is unique up to homothety. Each $m_x$ defines a measure on $\mathbb S^2$
which is invariant under
the maximal compact subgroup $\op{Stab}_G(x)$.

\begin{Def}[The Burger-Roblin measure]\label{BR} {\rm
  The Burger-Roblin measure $m^{\BR}_\G$ (\cite{Burger1990}, \cite{Roblin2003})
  associated to $\{\nu_x\}$ and $\{m_x\}$
is the measure on $\T^1(\G\ba \mathbb B)$
induced by the following $\G$-invariant measure on $\T^1(\mathbb B)$:
$$d  \tilde m^{\BR}(u)=e^{2\beta_{u^+}(x, \pi(u))}\;
 e^{\delta_\G \beta_{u^-}(x,\pi(u)) }\;
dm_x(u^+) d\nu_x(u^-) dt $$
for $x\in \mathbb B$.
%We denote by $m^{\BR}$ its induced measure on $\T^1(X)$ and call it
%the Burger-Roblin measure after those who first studied this measure \cite{Burger1990} and \cite{Roblin2003}.
By the conformal properties of $\{\nu_x\}$ and $\{m_x\}$,
this definition is independent of the choice of $x\in \mathbb B$.}
\end{Def}

%A limit point $\xi\in \Lambda(\G)$ is called a parabolic fixed point if
%it is stabilized by a parabolic element of $\G$.
 %We denote by $\Lambda_p(\G)$ the set of all parabolic fixed points
%of $\G$.

On $H/M=C_0^\dagger$, we consider the following two measures:
\begin{equation}\label{hars}d\mu^{\Leb}_{C_0^\dag} (s)=e^{2\beta_{s^+}(x,\pi(s))}dm_x(s)\quad\text{and}
\quad
 d\mu^{\PS}_{C_0^\dag} (s):=e^{\delta_\G \beta_{s^+}(x,\pi(s))} d\nu_x(s^+) \end{equation}
for $x\in \bB$.
These definitions are independent of the choice of $x$ and
$\mu^{\Leb}_{C_0^\dag}$ (resp. $\mu^{\PS}_{C_0^\dag}$) is left-invariant by $H$ (resp. $H\cap \G)$).
Hence we may consider the measures
$\mu^{\Leb}_{C_0^\dag}$ and $\mu^{\PS}_{C_0^\dag}$  on the quotient $(H\cap \G)\ba C_0^\dag$.

For the following theorem \ref{os} and the proposition \ref{ost},
we assume that the natural projection map $\G\cap H \ba \hat C_0 \to \G\ba \bB$ is proper and that
  any parabolic element of $\Gamma$
  fixing a point on $C_0$ stabilizes $C_0$.

\begin{Thm}\cite{OhShahGFH} \label{os}
  %$\Lambda_p(\G)\cap C_0\subset \Lambda_p(H\cap \G)$.
For $\psi\in C_c(\G\ba G/M)$,
as $t\to  \infty$,
$$e^{(2-\delta_\G)t}\int_{s\in (\G\cap H)\ba C_0^\dagger}\psi (sa_t) d\mu^{\Leb}_{C_0^\dag} (s)
\sim\frac{|\mu^{\PS}_{C_0^\dagger}|}{|m^{\BMS}_\G|}
m^{\BR}_\G(\psi)  $$
where $|\mu^{\PS}_{C_0^\dagger}|<\infty$.
Moreover $|\mu^{\PS}_{C_0^\dagger}| >0$ if $[\G: H\cap \G]=\infty$.
\end{Thm}

\begin{Prop}\cite{OhShahGFH} \label{ost} %  Suppose
  %that the natural projection map $\G\cap H \ba \hat C_0 \to \G\ba \bH^3$ is proper and that
  %$\Lambda_p(\G)\cap C_0\subset \Lambda_p(H\cap \G)$.
 For $\psi\in C_c(\G\ba G/M)$, there exists a compact subset
 $H_\psi\subset \G\cap H\ba H$ such that
 $$\psi(h a_t)=0$$
 for all $h\notin H_\psi$ and all $t\in \br$.
\end{Prop}

Letting $dm$ the probability invariant measure on $M$ and writing $h=sm \in C_0^\dag\times M$,
$dh=d\mu^{\Leb}_{C_0^\dag} (s) dm$ is a Haar measure on $H$, and
 the following defines a Haar measure on $G$: for $g=ha_tk\in HA^+K$,
$$dg=4\sinh 2t \cdot \cosh 2t  \; d h dt dm_{p_0}(k)$$
where $dm_{p_0}(k):=dm_{p_0}((k.X_0)^+)$.

We denote by $d\lambda$ the unique measure on $H\ba G$ which is compatible with the choice
of $dg$ and $dh$: for $\psi\in C_c(G)$,
$$\int_G\psi \; dg=\int_{[g]\in H\ba G}\int_{h\in H}\psi(h[g])\;dhd\lambda[g] .$$
Hence for $[e]=H$,
$d\lambda([e]a_tk)=4\sinh 2t \cdot \cosh 2t dtdm_{p_0}(k)$.

\section{Density of the Burger-Roblin measure on $\T^1_{p_0}(\bB)$ }
 Fixing $p_o, o\in \bB$, let $g_0\in G$ be such that $g_0(p_0)=o$.
 Let $\G<G$ be a non-elementary geometrically finite discrete subgroup of $G$.
 We use the same notation for
  $K=\op{Stab}_G(p_0)$, $A, A^+, X_0, M$, etc as in section \ref{not}.

  Let $N$ denote the expanding horospherical subgroup of $G$ for $A^+$:
$$N=\{g\in G: a_tg a_t^{-1}\to e\quad\text{as $t\to \infty$}\} .$$
The product map $A\times N\times K\to G$ is a diffeomorphism.

We fix a Borel subset $E\subset \S$ for the rest of this section.
\begin{Def}\rm
Define a function $\mathfrak R_E$ on $G$ as follows: for $g=a_tnk\in ANK$,
 $$\mathfrak R_E(g):=e^{ -\delta_\G t}\cdot \chi_{(g_0^{-1}E)_{p_0}}(k^{-1}) $$
where
$(g_0^{-1}E)_{p_0}:=\{u\in K: uX_0^-\in g_0^{-1}(E) \}$ and $\chi_{(g_0^{-1}E)_{p_0}}$
is its characteristic function. \end{Def}

\begin{Lem}\label{re} For any Borel subset $E\subset \S$,
$$\int_{k\in K/M} \mathfrak R_E(k^{-1} g_0)d\nu_{p_0}(k X_0^-)=\nu_o(E).$$
\end{Lem}
\begin{proof}
Write $k^{-1} g_0=a_t n k_0\in ANK$.
%Write $k^{-1} g_0=na_t k_0$ so that $\mathfrak R(k^{-1}g_0)=e^{-\delta_\G t}$.
Since $X_0^-=\lim_{s\to \infty}a_{-s}(p_0)$
and $\lim_{s\to \infty} a_{s+t} n a_{-s-t} = e$ ,
we obtain \begin{align*}& \beta_{kX_0^-}(o, p_0)
= \beta_{X_0^-}(k^{-1} o, p_0)\\ &=
\lim_{s\to \infty} d(k^{-1} g_0p_0, a_{-s}p_0) -d(p_0, a_{-s}p_0)
\\ &=
\lim_{s\to \infty} d(a_t n p_0, a_{-s}p_0) -d(p_0, a_{-s }p_0)
\\ &= \lim_{s\to \infty} d(( a_{s+t} n a_{-s-t}) a_{s+t} p_0, p_0) -d(p_0, a_{-s }p_0)
\\ &= \lim_{s\to \infty} ((s+t)-s) =t .\end{align*}
On the other hand, since $NA$ fixes $X_0^-$, $k_0^{-1} (X_0^-)=g_0^{-1} k (X_0^-)$,
 and hence $$ \chi_{(g_0^{-1}E)_{p_0}}(k^{-1}_0)=\chi_{g_0^{-1}E}(k_0^{-1}X_0^-)
=\chi_{E} (k(X_0^-)) .$$
 So
$$ \mathfrak R_E(k^{-1} g_0)=e^{-\delta_\G \beta_{k X_0^-}(o, p_0)} \chi_{E}(k(X_0^-)) .$$

 Therefore by the conformal property of
$\{\nu_x\}$,
\begin{align*}\int_{k\in K/M} \mathfrak R_E(k^{-1} g_0)
 d\nu_{p_0}(kX_0^-) &=
\int_{\xi\in  E}e^{-\delta_\G \beta_{\xi}(o, p_0)} d\nu_{p_0}(\xi)
=\nu_o (E).\end{align*}
\end{proof}

Fixing a left-invariant metric on $G$,
we denote by $U_\e$ an $\e$-ball around $e$, and for $S\subset G$, we set $S_\e=S\cap U_\e$.

\begin{Lem}( cf.  \cite[Lem. 6.1]{OhShahGFH}) \label{easy}
There exists $\ell\ge 1$ such
 that  for any $a_tnk\in ANK$ and small $\e>0$,
$$ a_tnk (g_0^{-1} U_{ \e} g_0) \subset A_{\ell \e} a_t N  K_{\ell \e}k . $$
\end{Lem}

For each small $\e>0$, we  choose a non-negative function $\psi^\e\in C_c(G)$ supported  inside
$U_{\e} $ and of integral $\int_{ G}\psi^\e dg$ one, and
define
$\Psi^\e\in C_c(\G\ba G)$ by
\begin{equation}\label{Pe}
\Psi^\e(g)=\sum_{\gamma\in \G}\psi^\e(\gamma g). \end{equation}
\begin{Def}\label{psiE} \rm Define a function $\Psi^\e_{E}$ on $\G\ba G$ by
$$\Psi^\e_{E}(g)=\int_{k^{-1} \in (g_0^{-1}E)_{p_0}}\Psi^\e(gkg_0^{-1})dm_{p_0}(k) .$$
\end{Def}
For each $\e>0$, define
\begin{equation}\label{ee}
  E^+_\e:=  g_0 U_\e g_0^{-1}(E)\quad\text{and}\quad E^-_\e :=\cap_{u\in U_\e} g_0 u g_0^{-1}(E) .\end{equation}

By Lemma \ref{easy}, it follows that there exists $c>0$ such that
for all $g\in U_\e$ and $g_1\in G$,
\begin{equation}\label{rrr} (1-c\e)
\mathfrak R_{E_\e^-}(g_1g_0)\le \mathfrak R_{E}(g_1 g g_0)\le (1+c\e)
\mathfrak R_{E_\e^+}(g_1g_0) .\end{equation}

\begin{Prop} \label{brd} There exists $c>0$ such that for all small $\e>0$,
$$ (1-c\e)  \nu_o(E^-_\e) \le m^{\BR}_\G(\Psi_E^\e)\le  (1+c\e) \nu_o (E^+_\e) .$$
\end{Prop}
\begin{proof}

Using the decomposition $G=ANK$,
we  have for $g=a_tnk$,
$$dg=dtdndm_{p_0}(k) $$
where $dn$ is the Lebesgue measure on $N$.

We use the following formula for $\tilde m^{\BR}$ (cf. \cite{OhShahGFH}): for any $\psi\in C_c( G)^M$,
$$\tilde m^{\BR}(\psi)=\int_{K}\int_{A}\int_N\psi( ka_t n)
e^{-\delta_\G t} \; dn dt d\nu_{p_0}(k(X_0^-)) .$$

 For $\psi^\e_E(g):=\int_{k^{-1} \in (g_0^{-1}E)_{p_0}}\psi^\e(gkg_0^{-1})dm_{p_0}(k) ,$
we have
\begin{align*}  &m^{\BR}_\G (\Psi_E^\e) =\tilde m^{\BR}(\psi_E^\e)\\&=
\int_{g\in G}\int_{k^{-1}\in (g_0^{-1}E)_{p_0}} \psi^\e(gk g_0^{-1}) dm_{p_0}(k)
d\tilde m^{\BR}(g)\\
&=\int_{KAN}\int_{k\in  K} \psi^\e (k_0 a_t n k g_0^{-1}) \chi_{(g_0^{-1}E)_{p_0}}(k^{-1}) dm_{p_0}(k)
e^{-\delta_\G t} dn dt d\nu_{p_0}(k_0X_0^-)
\\ &=\int_{k_0 \in K} \int_{ANK}  \psi^\e(k_0 (a_t n k) g_0^{-1}) \mathfrak R_E (a_t n k) dm_{p_0}(k)
 dn dt d\nu_{p_0}(k_0X_0^-)\\
&=\int_{k_0\in K} \int_{G} \psi^\e(k g g_0^{-1}) \mathfrak R_E (g) \, dg\;
 d\nu_{p_0}(k_0X_0^-)
\\
&=\int_{k_0\in K} \int_{G} \psi^\e( g ) \mathfrak R_E(k^{-1} g g_0)
dg\;  d\nu_{p_0}(k_0X_0^-). \end{align*}

Hence applying \eqref{rrr}, the identity $\int\psi^\e dg=1$ and Lemma \ref{re}, we deduce that
\begin{align*}  m^{\BR}_\G (\Psi_E^\e) &\le (1+ c\e)  \int_{k\in K} \left(\int_{G} \psi^\e( g ) dg \right)\mathfrak R_{E^+_\e}(k^{-1} g_0)
 d\nu_{p_0}(kX_0^-)  \\ &=
  (1+ c\e) \int_{k\in K} \mathfrak R_{E^+_\e}(k^{-1} g_0)
 d\nu_{p_0}(kX_0^-)
\\ &=  (1+ c\e) \nu_o(E^+_\e)  .\end{align*}
The other inequality follows similarly.
 \end{proof}

\section{Simpler proof of Theorem \ref{t_main} for the special case of $E=\S$.}
The result in this section is covered by the proof of
Theorem \ref{t_main} (for general $E$)
given in section \ref{finals}. However we present a separate proof
for this special case as
 it is considerably simpler and uses a different
interpretation of the counting function.

We may assume without loss of generality that $\P=\G(C_0)$. We use the notations from section \ref{not}.

Set
$$\N_T({\mathcal P})=\#\{C\in \P: d(\hat C, o)<t\}.$$
\begin{Lem} For $T>1$,
$$\N_T({\mathcal P})=\# [e]\G \cap [e] A_T^+K g_0^{-1}$$
where $[e]=H\in H\ba G$ and $A_T^+=\{a_t: 0\le t \le  T\}$.
\end{Lem}
\begin{proof}
%This observation is basically made in \cite[Pf of Thm. 2.5]{EskinMcMullen1993}.
Note that $N_T({\mathcal P})$ is equal to the number
of hyperbolic planes $\gamma (\hat C_0)$ such that $d(o, \gamma(\hat C_0))<T$,
or equivalently, $d(\gamma^{-1}(o),\hat C_0)<T$.
Since $\{x\in \mathbb B: d(x,\hat C_0)<T\}=HA_T^+(p_0)$,
 $N_T({\mathcal P})$ is equal to
 the number of $[\gamma]\in \G/\op{Stab}_{\G}({\hat C_0})$ such that
 $\gamma^{-1} g_0 p_0
 \in HA_T^+p_0$, or
 alternatively,
 the number of $[\gamma]\in H\cap \G\ba \G$ such that
 $\gamma g_0\in HA_T^+ K$, which is equal to
$ \# [e]\G g_0\cap [e] A_T^+K$.
\end{proof}

Define the following counting function $F_T$ on $\G\ba G$ by
$$F_T(g):=\sum_{\gamma\in \G\cap H\ba \G} \chi_{B_T}([e]\gamma g )$$
where $B_T=[e]A_T^+K g_0^{-1}\subset H\ba G$. Note that $F_T(e)=\N_T(\P)$.

By the strong wave front lemma (see \cite{GorodnikShahOhIsrael}),
for all small $\e>0$, there exists $\ell>1$ and $t_0>0$
such that for all $t>t_0$,
\be \label{stkak}
Ka_t k g_0^{-1} U_\e \subset Ka_t A_{\ell \e} k K_{\ell \e}  g_0^{-1} .\ee

%For all small $\e>0$, there exists $\ell>1$  (see \cite{EskinMcMullen1993},
% \cite{GorodnikShahOhIsrael})
It follows that for all $T\gg 1$,
\begin{equation*}\label{bu} (B_T -B_{t_0}) U_\e \subset B_{T+\ell \e}\;\;\text{and}\;\;
B_{T-\ell \e}-B_{t_0} \subset \cap_{u\in U_\e} B_Tu .\end{equation*}
Hence there exists $m_0 \ge 1$ such that for all $g\in U_\e$ and $T\gg 1$,
$$ F_{T-\ell \e}(g)-m_0\le   F_T(e)\le
 F_{T+\ell \e}(g)+m_0 .$$
Integrating against $\Psi^\e$ (see \eqref{Pe}), we obtain
\begin{equation}\label{e_in}\la F_{T-\ell\e},\Psi^\e \ra -m_0 \le F_T(e)\le
\la F_{T+\ell\e },\Psi^\e \ra +m_0 ,\end{equation}
where the inner product is taken with respect $dg$.

%\begin{Lem}\label{l_inner} For any $\psi\in C_c(\G\ba G)$,
%We have
%$$\la F_T, \psi\ra =\int_{k\in K}\int_{0\le t \le T} \int_{H\cap \G\ba C_0^\dag}
%\psi_{kg_0^{-1}}(sa_t)  (4 \sinh 2t \cdot \cosh 2t)  \;dsdrdm_{p_0}(k)$$
Setting $\Xi_t= 4\sinh 2t \cdot \cosh 2t$,
we have \begin{align*} &
\la F_{T+\ell \e}, \Psi^\e\ra \notag = \int_{g\in \G\cap H\ba G} \chi_{B_T} ([e] g) \Psi^\e(g g_0^{-1}) \; dg \notag \\
&=
\int_{k\in K} \int_{0}^{ T+\ell \e}
\int_{s\in \G\cap H\ba  C_0^\dag} \left(\int_{m\in M}
\Psi^\e(s  a_tm kg_0^{-1})\; dm\right) \Xi_t \;d\mu_{C_0^\dag}^{\Leb} (s) dtdm_{p_0}(k)
\notag \\ &=
\int_{k\in K} \int_{0}^{T+\ell \e }   \int_{s\in \G\cap H\ba C_0^\dag}
 \Psi_{kg_0^{-1}}^\e(sa_t) \Xi_t\;
d\mu_{C_0^\dag}^{\Leb} (s) dt dm_{p_0}(k). \notag
 \end{align*}
where $\Psi_{g_1}^\e\in C_c(\G\ba G)^M$ is given by
$\Psi_{g_1}^\e(g)=\int_{m\in M} \Psi^\e(gmg_1)\; dm $.

Hence by Proposition \ref{all} and Theorem \ref{os}, and using
$\Xi_t\sim e^{2t}$,
we deduce that as $T\to \infty$,
%(cf. Proof of prop. 6.8 in \cite{OhShahGFH})
\begin{equation*}\label{fin} \la F_{T+\ell\e}, \Psi^\e\ra
\sim \frac{  |\mu^{\PS}_{C_0^\dagger}|}{\delta_\G\cdot  |m^{\BMS}_\G|}
m^{\BR}_\G(\Psi^\e_{\S}) e^{\delta_\G (T+\ell \e)}\end{equation*}
where $$\Psi^\e_{\S}(g)=\int_{k\in K}\Psi^\e(gkg_0^{-1})dm_{p_0}(k) .$$

By Prop. \ref{brd},
$$m^{\BR}_\G(\Psi^\e_{\S})=(1+O(\e)) |\nu_o| .$$
Therefore it follows, as $\e>0$ is arbitrary,
$$\limsup_T\frac{F_T(e)}{e^{\delta_\G T}}\le
\frac{|\nu_o|\cdot
|\mu^{\PS}_{C_0^\dagger}|}{\delta_\G\cdot  |m^{\BMS}_\G|} .$$
Similarly
$$\liminf_T\frac{F_T(e)}{e^{\delta_\G T}}\ge
\frac{|\nu_o|\cdot
|\mu^{\PS}_{C_0^\dagger}|}{\delta_\G\cdot  |m^{\BMS}_\G|} .$$
This
finishes the proof, as $|\mu^{\PS}_{C_0^\dagger}|=\op{sk}_\G(\P)$.
\section{Uniform distribution along $\mathfrak b_T(\mathcal W)$}
 In this section, fix a Borel subset $\mathcal W\subset K$ with $M \W=\W$.

 \begin{Def}\label{btw}\rm  For $T>1$,
we set $$\mathfrak b_T(\mathcal W)= H\ba H KA^+_T \W \subset H\ba G $$
where $A_T^+=\{a_t\in A: 0\le t\le T\}$.
\end{Def}

%The goal of this section is to deduce the following theorem \ref{mtt} from Theorem \ref{os}.
%Most of technical difficulties lie in understanding the shape of the set $B_T(E)$ in $HA^+K$ decomposition.
%and how it changes when multiplied by small neighborhoods of $e$ of $G$.

We assume that
 the natural projection map $\G\cap H \ba \hat C_0 \to \G\ba \bB$ is proper and that
  any parabolic element of $\Gamma$
  fixing a point on $C_0$ stabilizes $C_0$.

\begin{Thm}\label{mtt}
For any $\psi\in C_c(\G\ba G)$, we have
$$\int_{g\in \mathfrak b_T(\W)}\int_{h\in \G\cap H\ba H}\psi(hg)dhd\lambda(g)
\sim \frac{|\mu^{\PS}_{C_0^\dagger}|}{\delta_\G \cdot |m^{\op{BMS}}|}\cdot
 \int_{k\in \W} m^{\BR}_\G(\psi_{k }  )\; dm_{p_0}(k)  \cdot e^{\delta_\G T}$$
as $T\to \infty$, where $\psi_k\in C_c(\G\ba G)^M$ is given by
$\psi_k(g)=\int_{m\in M}\psi(gmk)dm$.
\end{Thm}
\begin{proof} (cf. \cite[Thm 4.3]{OhShahcircle})

%Let $\ell=\ell(E)\ge 1$ be as in Proposition \ref{uk}.
Set $K_\e'=\cup_{k\in K} k K_\e k^{-1}$ and
 define $\psi_\e^{\pm}\in C_c(\G\ba G)$ by
$$\psi_\e^+(g):=\sup_{u\in K_\e' } \psi(gu)\quad \text{and }\quad \psi_\e^-(g):=\inf_{u\in K_\e' }\psi(gu).$$
Note that for a given $\eta>0$,
there exists $\e=\e(\eta)>0$ such that for all $g\in \G\ba G$,
$|\psi_\e^+(g)-\psi_\e^-(g)|\le \eta$ by the uniform continuity of $\psi$.

We can deduce from Theorem \ref{os} that
for all $t>T_1(\eta)\gg 1$,
$$\int_{h\in \G\cap H\ba H}\psi_\e^+(ha_t k ) dh =
(1+O(\eta)) \frac{|\mu^{\PS}_{C_0^\dagger}|}{|m^{\op{BMS}}_\G|} m^{\op{BR}}_\G(\psi_{\e,k}^+) e^{(\delta-2) t}$$
 where $\psi_{\e,k}^+$ is defined similarly as $\psi_k$
and the implied constant can be taken uniformly over all $k\in K$. %Similarly, for all $t<-T_1(\eta)$,
 %\begin{align}\label{pstars}&\int_{h\in \G\cap H\ba H}\psi_\e^+(ha_tn) dh =
%(1+O(\eta)) \frac{|\mu^{\PS}_{C_0^\ddagger}|}{|m^{\op{BMS}}|} m^{\op{BR}}(\psi_{\e,n}^+)
% e^{(\delta-2)|t|}\end{align*}
Defining $$K_T(t):=\{k\in K: a_t k\in HKA_T^+\} ,$$
by Prop. 4.8 and Corollary 4.11 in \cite{OhShahcircle},
we have $H KA_T^+=\cup_{0\le t\le T} Ha_t K_T(t) $
and
there exists a sufficiently large $T_0(\e)>T_1(\eta)$ such that
$e\in K_T(t)\subset K_\e M$ for all $T_0(\e)<t<T$.

%Let $T_0(\e)>T_1(\eta)$ be as in
% Proposition \ref{ksmall}.
For $[e]=H\in H\ba G$ and $s>0$,
set $$V_T(s):=(\cup_{s\le t\le  T}[e]a_tK_T(t))\W$$
so that $$\mathfrak b_T(\W )=V_T(s)\cup (\mathfrak b_T(\W)- V_T(s)).$$

 %$$V_T(T_0(\e)):=V_T(T_0(\e))$$ and
%Setting $$\psi^H(g):= \int_{h\in \G\cap H\ba H}\psi(hg) dh ,$$

Let $[g]=[e]a_tkk_1\in V_T(T_0(\e))$ where $k_1\in K$ and $k\in \W$.
For $t>T_0(\e)$,
there exist $h_0\in H$ and $ u\in K_{\e }'$ such that $a_tk_1 k =h_0a_t k u$
and hence
\begin{multline*}
\psi^H(g):=\int_{h\in \G\cap H\ba H}\psi(h g) dh \\
=\int_{h\in \G\cap H\ba H}\psi(hh_0 a_t k u) dh \le
\int_{h\in \G\cap H\ba H}\psi_\e^+(ha_tk) dh. \end{multline*}

Therefore
$$\int_{V_T(T_0(\e))}\psi^H(g) d\lambda(g) \le
\int_{k\in \W} \int_{T_0(\e)<t \le  T} \int_{h\in \G\cap H\ba H}\psi_\e^+(ha_t k)
\Xi_t  \; dh dtdm_{p_0}(k) $$
where $\Xi_t=4\sinh 2t \cosh 2t$.

Using $\Xi_t\sim e^{2t}$, we then deduce
\begin{multline*}
\int_{k\in \W}\int_{T_0(\e) < t< T}\int_{h\in \G\cap H\ba H}\psi_\e^+(ha_tk) \Xi_t\; dh  dt dm_{p_0}( k)
\\ = (1+O(\eta)) \frac{|\mu^{\PS}_{C_0^\dagger}|}{\delta_\G \cdot |m^{\op{BMS}}|}\cdot
 \int_{k\in \W} m^{\BR}_\G(\psi_k ) dm_{p_0}( k)  \cdot
  ( e^{\delta_\G T} -e^{\delta_\G T_0(\e)})
\end{multline*}
since  $m^{\BR}_\G(\psi_{\e, k}^+)=(1+O(\eta))m^{\BR}_\G(\psi_k)$.
%  where the implied constant depends only on $\psi$.

Hence
\begin{multline*}
\limsup_T  \frac{ \int_{ V_T(T_0(\e))} \psi^H(g) d\lambda(g)}{e^{\delta_\G T}}
=
(1+O(\eta)) \frac{|\mu^{\PS}_{C_0^\dagger}|}{\delta_\G \cdot |m^{\op{BMS}}|}\cdot
 \int_{\W} m^{\BR}_\G(\psi_k ) dm_{p_0}( k) .\end{multline*}

On the other hand, it follows from Proposition \ref{ost}
that $$ \int_{[g]\in \mathfrak b_T(\W)-
V_T(T_0(\e))} \int_{h\in \G\cap H\ba H}\psi(hg) dh d\lambda(g) =O(1) .$$
As $\eta>0$ is arbitrary and $\e(\eta)\to 0$ as $\eta\to 0$, it follows that
$$\limsup_T
\frac{\int_{[g]\in \mathfrak b_T(\W)} \psi^H(g) d\lambda(g)}
{e^{\delta_\G T}}\le  \frac{|\mu^{\PS}_{C_0^\dagger}|}{\delta_\G \cdot |m^{\op{BMS}}|}\cdot
 \int_{k\in \W} m^{\BR}_\G(\psi_k ) dm_{p_0}(k) .$$
By a similar argument, one can prove
$$\liminf_T
\frac{\int_{[g]\in \mathfrak b_T(\W)} \psi^H(g) d\lambda(g)}{e^{\delta_\G T}}
\ge  \frac{|\mu^{\PS}_{C_0^\dagger}|}{\delta_\G \cdot |m^{\op{BMS}}|}\cdot
 \int_{k\in \W} m^{\BR}_\G(\psi_k ) dm_{p_0}(k) .$$
\end{proof}

\section{Proof of Theorem \ref{t_main}}\label{finals}
Without loss of generality, we may assume that $\P=\G(C_0)$.
 We keep the notations from section \ref{not}.

%We set $K_0=g_0Kg_0^{-1}, M_0=g_0Mg_0^{-1}, A_0^{\pm}=g_0A^{\pm} g_0^{-1}, H_0=g_0^{-1}H g_0$
%and note that $g_0(X_0)$ is a vector based at $o$.
%Then $T_{o}(\mathbb B)=K_0/M_0$

%In the following we assume that $\G(C_0)$ is a locally finite circle packing of $\mathbb S^2$.
\begin{Def}\label{adm}{\rm  A subset $E\subset \S$ is said to be $\mathcal P$-admissible if,
  for any $C\in \mathcal P$, $C^\circ\cap E\ne \emptyset$ implies
  $C^\circ\subset E$, possibly except for finitely many circles.
}\end{Def}

For a subset $E\subset \S$, we define
$E_{p_0}\subset K$ by
$$E_{p_0}:=\{k\in K: k(X_0^-)\in E\} .$$
We also set
$$\mathcal N_T(\P, E):=\{C\in \P: C\cap E\ne \emptyset,\;\; d(\hat C, o)<T\} .$$

\begin{Lem}\label{l_ad} Fix a $\P$-admissible subset $E\subset \S$.
 For all $T>1$,
$$\N_T({\mathcal P},E)=\# [e]\G  \cap [e]K A_T^+ (g_0^{-1}E)_{p_0}^{-1} g_0^{-1} $$
up to a uniform finite additive constant where $[e]=H\in H\ba G$.
% and $g_0\in G$ such that $o=g_0(p_0)$.
\end{Lem}
\begin{proof} Since $g_0K/M$ represents the set of all unit vectors based at $o$,
and the set $\{u\in \T_o^1(\mathbb B): u^-\in E\}$
is identified with $g_0 (g_0^{-1}E)_{p_0}=\{g_0k[M]: kX_0^-\in g_0^{-1}E\}$,
the set $g_0 (g_0^{-1}E)_{p_0}A^- (p_0)$
 represents the set of all points in $\mathbb B$
lying in the cone consisting of geodesic rays connecting $o$ with a point in $E$.
Therefore the condition $C\subset E$ is equivalent to that $\hat C\subset g_0 (g_0^{-1}E)_{p_0}A^- (p_0)$.
Hence by the $\P$-admissibility condition,
we may assume without loss of generality that $\N_T({\mathcal P},E)$ is equal to the number
of hyperbolic planes $\gamma (\hat C_0)$ such that $d(o, \gamma(\hat C_0))<T$
and $\gamma (\hat C_0)\subset  g_0 (g_0^{-1}E)_{p_0}A^- (p_0)$.
Since $\{x\in \mathbb B: d(o,x)<T\}= g_0KA_T^-(p_0)$
where $A_T^-=\{a_{-t}:0\le t\le T\}$,
the former condition is again same as $\gamma(\hat C_0)\cap g_0 KA_T^-(p_0)\ne\emptyset$.
%If we identify $\B^3$ with $G/K_0$,
%we have $\hat C_0=g_0^{-1} H_0/H_0\cap K_0$,
%and
Hence
 \begin{align*}
 \N_T({\mathcal P},E)&=\{\gamma(C_0): \gamma(\hat C_0)\cap g_0 KA_T^-(p_0)\ne\emptyset,\;
 \gamma (\hat C_0) \subset g_0 (g_0^{-1}E )_{p_0}A^- (p_0) \}
\\
&=\{ [\gamma] \in \G/ \Gamma\cap H: \gamma \in  g_0 KA_T^-K H\cap  g_0 (g_0^{-1}E)_{p_0}A^- KH\}\\
&=\{ [\gamma] \in \G/ \Gamma\cap H: \gamma \in   g_0 (g_0^{-1}E)_{p_0}A^-_T KH\}.\end{align*}
In the last equality, we have used the fact that
if $a_{-t}\in KA_T^-KH$ for some $t>0$,  then $t<T$ (see
\cite[Lem 4.10]{OhShahcircle}).

By taking the inverse, we obtain that
$$\N_T({\mathcal P},E) =[e]\G\cap [e] K A_T^+ (g_0^{-1} E)_{p_0}^{-1} g_0^{-1} .$$
\end{proof}

Fixing a Borel subset $E\subset \S$,
recall the definition of $E_{\e}^{\pm}$ from \eqref{ee}:
\begin{equation*}
  E^+_\e:=  g_0 U_\e g_0^{-1}(E)\quad\text{and}\quad E^-_\e :=\cap_{u\in U_\e} g_0 u g_0^{-1}(E) .\end{equation*}

We can find a $\P$-admissible Borel subset $\tilde E_{\e}^{+}$ such that
$E\subset \tilde E_{\e}^{+}\subset E_{\e}^+$  by adding
all the open disks inside $E_{\e}^{+}$
 intersecting the boundary of $E$.
  Similarly we can find a $\P$-admissible Borel subset $\tilde E_{\e}^{-}$ such that
 $E_{\e}^-\subset\tilde E_{\e}^-\subset E$ by adding all the open disks inside $E$
 intersecting the boundary of $E_{\e}^{-}$.
 By the local finiteness of $\P$, there are only finitely many circles intersecting
  $\partial(E)$ (resp.
 $\tilde E_{\e}^{-}$) which are not contained in $\tilde E_{\e}^{+}$ (resp. $E$). Therefore
there exists $q_\e\ge 1$ (independent of $T$) such that
\begin{equation}\label{npt}\N_T(\P, \tilde E_{\e}^{-} ) - q_\e\le \N_T(\P, E) \le \N_T(\P, \tilde E_{\e}^{+} ) +q_\e .\end{equation}

%Define the following counting function $F_T$ on $\G\ba G$ by
%$$F_T(g):=\sum_{\gamma\in \G\cap H\ba \G} \chi_{B_T}([e]\gamma g )$$
%where
% $$B_T(E):= [e]KA^+_T(g_0^{-1} E)_{p_0}^{-1} g_0^{-1} \subset H\ba G .$$

%Note that $F_T(e)=\N_T(\P, E)$.

 Setting
 $$B_T(E):= [e]KA^+_T(g_0^{-1} E)_{p_0}^{-1} g_0^{-1} \subset H\ba G ,$$
we define functions $F_T^{\e, \pm}$ on $\G\ba G$:
 $$F_T^{\e, \pm}(g):=\sum_{\gamma\in \G\cap H\ba \G} \chi_{B_{T\pm \ell \e} (E_{(\ell +1)\e}^{\pm})}([e]\gamma g ) .$$

\begin{Lem} There exists $m_\e\ge 1$ such that  for all $g\in U_{ \e }$ and $T\gg 1$,
\begin{equation*}\label{ft} F_T^{\e, +}(g)-m_\e
\le  \N_T(\P, E) \le  F_T^{\e, +}(g)+m_\e.
\end{equation*}\end{Lem}
\begin{proof}
It follows from \eqref{stkak} that
$$B_T(E_\e^+) U_\e \subset B_{T+\ell \e} (E_{(\ell+1) \e}^{+})\;\;\text{and}\;\;
B_{T-\ell \e} (E_{(\ell +1)\e}^{-})\subset \cap_{u\in U_\e} B_T(E_\e^-) u .$$
Hence for any $g\in  U_{\e }$, as $U_\e$ is symmetric,
$$\# [e]\G\cap B_T(\tilde E^+_\e) \le \# [e]\G\cap B_T(\tilde E^+_\e) U_{ \e} g^{-1} \le
\# [e]\G g \cap B_{T+\ell \e}(E_{(\ell +1)\e}^+) .$$
By Lemma \ref{l_ad} and \eqref{npt}, it follows that for some fixed $m_\e \ge 1$,
$$ \N_T(\P, E) \le  F_T^{\e, +}(g)+m_\e.$$
 The other inequality can be proved similarly.
\end{proof}

 Hence by integrating against $\Psi^\e$ (see \eqref{Pe}), we obtain
\begin{equation}\label{e_in}\la F_{T}^{\e, -},\Psi^\e \ra - m_\e \le \N_T(\P, E)\le
\la F_{T}^{\e, +},\Psi^\e \ra +m_\e.\end{equation}

We note that
$$B_T(E)=\mathfrak b_T( (g_0^{-1}E)_{p_0}^{-1})\; g_0^{-1}$$
where $\mathfrak b_T(\W)$ is defined as in Def. \ref{btw}.

Since \begin{align*}
\la F_{T}^{\e,+}, \Psi^\e \ra \notag  &=\int_{\G\cap H \ba G}
 \chi_{B_{T+\ell \e}(E_{(\ell +1)\e}^+)}([e] g )
\Psi^\e(g) \; dg\\
\\ &=\int_{[g]\in B_{T+\ell \e}(E_{(\ell +1)\e}^+)}\int_{h\in \G\cap H\ba H}\Psi^\e (hg)dhd\lambda(g)
\\ &=\int_{[g]\in \mathfrak b_{T+\ell \e}( (g_0^{-1}E_{(\ell +1)\e}^+)_{p_0}^{-1})}\int_{h\in \G\cap H\ba H}\Psi^\e (hg g_0^{-1})dhd\lambda(g)
\end{align*}
we deduce from Proposition \ref{all} and Theorem \ref{mtt}
that
\begin{equation}\label{sim2} \la F_{T}^{\e,+}, \Psi^\e \ra \sim
\frac{\op{sk}_{\G}(C_0)}{\delta_\G \cdot |m^{\op{BMS}}_\G|}\cdot
  m^{\BR}_\G(\Psi_{E_{(\ell +1)\e}^+}^\e )  \cdot e^{ \delta_\G (T+\ell \e)} \end{equation}
where $\Psi^\e_{E}(g)=\int_{k^{-1} \in (g_0^{-1}E)_{p_0}}\Psi^\e(gkg_0^{-1})dm_{p_0}(k) $
 (see Def. \ref{psiE}).

Therefore by \eqref{e_in} and Prop. \ref{brd}
we have
$$\limsup_T \frac{ \N_T(\P, E)}{ e^{\delta_\G T}}\le
(1+O(\e)) \frac{\op{sk}_{\G}(C_0)}{\delta_\G \cdot |m^{\op{BMS}}_\G|}\cdot \nu_o(E_{(\ell+1)\e}) .$$
Since $\nu_o(\partial(E))=0$ by the assumption,
$\nu_o(E_{(\ell+1)\e} -E)\to 0$ as $\e \to 0$.
As $\e$ can be taken arbitrarily small,
it follows that
$$\limsup_T \frac{ \N_T(\P, E)}{ e^{\delta_\G T}}\le
 \frac{\op{sk}_{\G}(C_0)}{\delta_\G \cdot |m^{\op{BMS}}_\G|}\cdot \nu_o(E) .$$
Similarly, we can prove
$$\liminf_T \frac{ \N_T(\P, E)}{ e^{\delta_\G T}}\ge
 \frac{\op{sk}_{\G}(C_0)}{\delta_\G \cdot |m^{\op{BMS}}_\G|}\cdot \nu_o(E) .$$
This completes the proof.

\bibliographystyle{plain}
%\bibliography{Ohbibliog}

\begin{thebibliography}{10}

\bibitem{Ahlfors1964}
Lars~V. Ahlfors.
\newblock Finitely generated {K}leinian groups.
\newblock {\em Amer. J. Math.}, 86:413--429, 1964.

\bibitem{Bowen1971}
Rufus Bowen.
\newblock Periodic points and measures for {A}xiom {$A$} diffeomorphisms.
\newblock {\em Trans. Amer. Math. Soc.}, 154:377--397, 1971.

\bibitem{Burger1990}
Marc Burger.
\newblock Horocycle flow on geometrically finite surfaces.
\newblock {\em Duke Math. J.}, 61(3):779--803, 1990.

\bibitem{Claytor1934}
Schieffelin Claytor.
\newblock Topological immersion of {P}eanian continua in a spherical surface.
\newblock {\em Ann. of Math. (2)}, 35(4):809--835, 1934.

\bibitem{ErickssonLagarias2007}
Nicholas Eriksson and Jeffrey~C. Lagarias.
\newblock Apollonian circle packings: number theory. {II}. {S}pherical and
  hyperbolic packings.
\newblock {\em Ramanujan J.}, 14(3):437--469, 2007.

\bibitem{FlaminioSpatzier}
L.~Flaminio and R.~J. Spatzier.
\newblock Geometrically finite groups, {P}atterson-{S}ullivan measures and
  {R}atner's rigidity theorem.
\newblock {\em Invent. Math.}, 99(3):601--626, 1990.

\bibitem{GorodnikShahOhIsrael}
Alex Gorodnik, Nimish Shah, and Hee Oh.
\newblock Strong wavefront lemma and counting lattice points in sectors.
\newblock {\em Israel J. Math}, 176:419--444, 2010.

\bibitem{Kojima1992}
Sadayoshi Kojima.
\newblock Polyhedral decomposition of hyperbolic {$3$}-manifolds with totally
  geodesic boundary.
\newblock In {\em Aspects of low-dimensional manifolds}, volume~20 of {\em Adv.
  Stud. Pure Math.}, pages 93--112. Kinokuniya, Tokyo, 1992.

\bibitem{KontorovichOh}
Alex Kontorovich and Hee Oh.
\newblock Apollonian circle packings and closed horospheres on hyperbolic
  3-manifolds.
\newblock {\em Preprint, 2009}.

\bibitem{LagariasMallowsWilks}
Jeffrey~C. Lagarias, Colin~L. Mallows, and Allan~R. Wilks.
\newblock Beyond the {D}escartes circle theorem.
\newblock {\em Amer. Math. Monthly}, 109(4):338--361, 2002.

\bibitem{MardenOutercircles}
A.~Marden.
\newblock {\em Outer circles}.
\newblock Cambridge University Press, Cambridge, 2007.
\newblock An introduction to hyperbolic 3-manifolds.

\bibitem{Margulisthesis}
Gregory Margulis.
\newblock {\em On some aspects of the theory of {A}nosov systems}.
\newblock Springer Monographs in Mathematics. Springer-Verlag, Berlin, 2004.
\newblock With a survey by Richard Sharp: Periodic orbits of hyperbolic flows,
  Translated from the Russian by Valentina Vladimirovna Szulikowska.

\bibitem{MumfordSeriesWright}
David Mumford, Caroline Series, and David Wright.
\newblock {\em Indra's pearls}.
\newblock Cambridge University Press, New York, 2002.
\newblock The vision of Felix Klein.

\bibitem{OhICM}
Hee Oh.
\newblock Dynamics on {G}eometrically finite hyperbolic manifolds with
  applications to {A}pollonian circle packings and beyond.
\newblock {\em To appear in the Proc. of I.C.M (Hyperabad, 2010)}.

\bibitem{OhShahcircle}
Hee Oh and Nimish Shah.
\newblock The asymptotic distribution of circles in the orbits of {K}leinian
  groups.
\newblock {\em Preprint}.

\bibitem{OhShahGFH}
Hee Oh and Nimish Shah.
\newblock Equidistribution and counting for orbits of geometrically finite
  hyperbolic groups.
\newblock {\em Preprint}.

\bibitem{Patterson1976}
S.J. Patterson.
\newblock The limit set of a {F}uchsian group.
\newblock {\em Acta Mathematica}, 136:241--273, 1976.

\bibitem{Roblin2003}
Thomas Roblin.
\newblock Ergodicit\'e et \'equidistribution en courbure n\'egative.
\newblock {\em M\'em. Soc. Math. Fr. (N.S.)}, (95):vi+96, 2003.

\bibitem{Schlichtkrull1984}
Henrik Schlichtkrull.
\newblock {\em Hyperfunctions and harmonic analysis on symmetric spaces},
  volume~49 of {\em Progress in Mathematics}.
\newblock Birkh\"auser Boston Inc., Boston, MA, 1984.

\bibitem{Sullivan1979}
Dennis Sullivan.
\newblock The density at infinity of a discrete group of hyperbolic motions.
\newblock {\em Inst. Hautes \'Etudes Sci. Publ. Math.}, (50):171--202, 1979.

\bibitem{Sullivan1984}
Dennis Sullivan.
\newblock Entropy, {H}ausdorff measures old and new, and limit sets of
  geometrically finite {K}leinian groups.
\newblock {\em Acta Math.}, 153(3-4):259--277, 1984.

\bibitem{Thurstonbook}
William Thurston.
\newblock {\em The Geometry and Topology of Three-Manifolds}.
\newblock available at www.msri.org/publications/books. Electronic
  version-March 2002.

\end{thebibliography}

\end{document}